\documentclass[12pt]{article}

\usepackage{amscd}
\usepackage{amssymb,amsmath}

\newtheorem{theorem}{Theorem}
\newtheorem{prop}[theorem]{Proposition}
\newtheorem{lemma}[theorem]{Lemma}
\newtheorem{cor}[theorem]{Corollary}

\author{Vladimir Baranovsky}
\title{Norm functors and effective zero cycles}
\date{January 12, 2009, revised February 23, 2009}

\begin{document}
\maketitle

\begin{abstract}
We compare two known definitions for a relative family of 
effective zero cycles, based
on traces and norms of functions, respectively. In characteristic zero we 
show that both definitions agree. In the general setting, we show 
that the norm map on functions can be expanded to a
norm functor between certain 
categories of line bundles, therefore giving a third 
approach to families of zero cycles.
\end{abstract}


\section{Introduction}

Let us start with a simple situation of a quasi-projective scheme $X$
over a perfect field $k$. We want to understand the notion of a family of 
degree $d$
effective zero cycles  parameterized by a $k$-scheme $T$.
When $T = Spec(k)$ these are just finite formal linear combinations
\begin{equation}
\label{cycle}
\xi = \sum d_i x_i
\end{equation}
of closed points $x_i \in X$ with $d_i \in \mathbb{Z}_{\geq 0}$. 
Such a cycle has degree $d = \sum_i d_i e_i$ where
$e_i$ is the degree of the field extension $k \subset k(x_i)$. 
We remark here that if $k \subset k(x_i)$ were not separable one would have 
to work with rational $d_i$ having powers of
$p = char\; k$ in denominators, cf. Section 8 of \cite{R2}.

To obtain a description which works over an arbitrary base $T$ let 
$f$ be a regular function defined on an open subset $U \subset X$ 
containing all $x_i$. Define the following elements in $k$
$$
\theta(f) = \sum_i d_i \theta_i (f(x_i)); \qquad n(f) = \prod_i 
n_i(f(x_i))^{d_i}
$$
where $\theta_i$ and $n_i$ stands for the trace and the norm of 
the field extension $k \subset k(x_i)$, respectively. 

When $\xi$ varies with $t \in T$, denote $X_T= X\times_k T$ and let
$\pi_T: X_T \to T$ be the canonical projection. The above construction 
gives \textit{trace} and \textit{norm} maps
$$
\theta: (\pi_T)_* \mathcal{O}_{\widehat{X}_T} \to \mathcal{O}_T; 
\qquad n: (\pi_T)_* \mathcal{O}_{\widehat{X}_T} \to \mathcal{O}_T.
$$
where $\widehat{X}_T$ is the completion of $X_T$ along the closed 
subset $Z$ swept out by the points $x_i(t)$. Observe
that $\theta$ is a morphism of $\mathcal{O}_T$-modules, while 
$n$ is just a multiplicative map. The values on functions pulled back from
 $T$ are given by $\theta(f) = d\cdot f$ and $n(f) = f^d$. 
Both trace and norm  should be continuous, i.e. 
factor through $(\pi_T)_* \mathcal{O}_Y$
for some closed subscheme $Y \hookrightarrow X$ with support in $X$.
Both constructions should commute with base change $T' \to T$. 

\bigskip
\noindent
This suggests the idea that a family of zero cycles with base $T$ can be 
defined by specifying an appropriate closed subset $Z \subset X_T$ 
and either a trace map $\theta$ or a norm map $n$, as above
(the version with norm maps is apparently originally due to 
Grothendieck, who also applied it to non-effective cycles by 
restricting to multiplicative groups of $\mathcal{O}_{\widehat{X}_T}$ and 
$\mathcal{O}_T$). Observe 
that for a base change $T' \to T$ the trace map automatically pulls back
to $T'$, being a morphism of $\mathcal{O}_T$-modules, while with the norm 
map $n$ we have to \textit{specify} the pullback of $n$. In other words, 
we should  have a system of maps 
$$
n_{T'}: (\pi_{T'})_* \mathcal{O}_{\widehat{X}_{T'}} \to \mathcal{O}_{T'}
$$
for $T' \to T$, which agree with each other in a natural sense. This seems
like an inconvenient detail. However, it turns out that one needs to 
impose further conditions on $\theta$ to get a trace map which comes 
from a geometric family of cycles (cf. definition 
after Lemma \ref{trace} in Section 3.1 below) while for $n$ the existence
of its extensions $n_{T'}$ plays the role of such a condition. 
In addition, the trace construction 
only works well when $k$ has characteristic zero (or finite characteristic
$p > d$).

The approach using traces was used in characteristic zero by B. Angeniol, 
cf. \cite{An}, and also by Buchstaber and Rees, cf. \cite{BR}. 
Angeniol extends his definition to cycles of higher dimensions, which 
leads to a construction of the Chow scheme of cycles. In the affine case
the norm approach was used by N. Roby in \cite{Ro}, but the global version 
of his construction was carried out only recently (more than 40 years later!) 
by D. Rydh in \cite{R1}-\cite{R3}. The latter author
deals with a general situation of a separated morphism of 
algebraic spaces $\pi: X \to S$. He also considers higher dimensional 
cycles, using an old idea of Barlet that a family of
$n$ dimensional cycles over an $l$-dimensional base can be represented locally
as a family of zero cycles over an $n+l$ dimensional base.

In the above construction one should take into account that 
$\theta$ and $n$ could factor through a completion along a smaller
closed subset $Z' \subset Z$. In the additive case, B. Angeniol formulates
a non-degeneracy condition ensuring that such $Z'$ does not exist.
In the multiplicative case, D. Rydh simply considers pairs $(Y, n)$ consisting
 of a closed subscheme $Y$ and a norm 
map $n: (\pi_T)_* \mathcal{O}_Y\to \mathcal{O}_T$, and then uses an equivalence
relation which identifies $(Y', n')$ and $(Y'', n'')$ if $n'$ and $n''$ 
factor through a norm map defined on a closed subscheme $Y\subset Y' \cap Y''$. 
In this paper we adopt the second approach, modifying and 
generalizing the trace definition.

In characterisic zero both functors of families of zero cycles are 
represented by the symmetric power $Sym^d (X/S)$, i.e. the quotient
of the $d$-fold cartesian product of $X$ over $S$ by the action 
of the symmetric group $\Sigma_d$. In arbitrary characteristic the 
approach based on traces breaks down: e.g.  we are not able to 
distinguish $\xi = x$ from $\xi = (p+1)x$, while the norms 
approach leads to the \textit{space of divided powers} $\Gamma^d(X/S)$.
There is a natural morphism $Sym^d(X/S) \to \Gamma^d(X/S)$ which is
an isomorphism in characteristic zero, but in general only a 
universal homeomorphism, cf. \cite{R1}. Even for 
general schemes over a field $k$ (not necessarily quasi-projective) 
both $Sym^d(X/S)$ and $\Gamma^d(X/S)$ may not be schemes but rather
algebraic spaces. Therefore, it is more natural to work in the category 
of algebraic spaces from the  beginning.

\bigskip
\noindent
The purpose of this paper is to explore a third approach to families of zero cycles, which
admits a reasonably straightforward generalization to higher dimensional
 cycles. In the original setup of a scheme $X$ over $k$, choose a line bundle $L$
defined on an open subset $U \subset X$ containing all $x_i$ 
in \eqref{cycle}
 and assume for simplicity that all $k(x_i)$ 
are equal to $k$. Define a one-dimensional vector space
$$
N(L) = \prod L_{x_i}^{\otimes d_i}
$$
over $k$. When $\xi$ varies over a base $T$, this gives a line bundle $N(L)$
on $T$. Obviously, an isomorphism of line bundles on $X$ induces an 
isomorphism of bundles on $T$. However, 
non-negativity of the coefficients  $d_i$ is reflected
in the fact that any \textit{any} morphism
of invertible $\mathcal{O}_{\widehat{X}_T}$-modules
$\psi: L \to M$   gives a morphism of $\mathcal{O}_T$-modules
$N(\psi): N(L) \to N(M)$. We can further consider the line bundles defined
only on a neighborhood of $Z$. Thus, for a closed subscheme 
$Y \hookrightarrow X_T$ supported at $Z$ we should have a
\textit{norm functor}
$$
N: PIC(Y) \to PIC(T)
$$
where $PIC$ is the category of line bundles and morphisms as 
$\mathcal{O}$-modules. Again, this
functor should come with functorial pullbacks with respect to morphisms 
 of schemes (or algebraic spaces) $T' \to T$. In practice it suffices to 
restrict to those $T'$ which are affine over $T$ (or even to the full 
subcategory generated by affine spaces over $T$).

However, on morphisms the correspondence
$\psi \mapsto N(\psi)$ is no longer $\mathcal{O}_T$-linear, but rather 
satisfies 
$N(\psi)(\pi_T^*(f) s) = f^d N(\psi)(s)$ where $s$ is a local section of
$N(L)$ and $f$ is a local section of $\mathcal{O}_T$. Morphisms of modules
with this property where also considered by N. Roby, cf. \cite{Ro} where
they are called \textit{homogeneous polynomial laws of degree $d$}. As with
norm maps, we should also specify a functorial extension of $\psi \mapsto N(\psi)$ with respect to base changes $T' \to T$. The fact that
$N$ is a functor means that $\psi \mapsto N(\psi)$ is multiplicative since
compositions should go to compositions. 
In addition, $N$ should agree with  tensor products of line bundles and, 
similarly to identities:
$$
\theta(\pi_T^*(f)) = d f; \qquad n(\pi_T^*(f)) = f^d; 
$$
we should have an isomorphism of functors
$$
\eta: N \circ \pi_T^* \simeq \{L \mapsto L^{\otimes d}\} 
$$
agreeing with base change. This rigidification also ensures that 
$N$ does not have any non-trivial functor automorphisms. 

Besides the generalization to higher dimensional cycles based on  the
work of F. Ducrot, R. Elkik and E. Mu\~noz-Garcia, cf. \cite{Du}, 
\cite{El}, \cite{MG} this approach to zero cycles also can be used to define the Uhlenbeck compactification of
moduli stack of vector bundles on a surface. The standard constructions
like Hilbert-to-Chow morphism, sums of cycles and Chow forms are
also rather simple in the language of norm functors. 

Norms of line bundles were earlier considered in \cite{EGA2} and 
\cite{De}. More general norms of quasi-coherent sheaves were
studied in \cite{Fe} and \cite{R2}.

\bigskip
\noindent
This work is organized as follows. In Section 2 we recall the basic results on
polynomial laws, divided powers and norms for finite flat morphisms. In 
Section 3 we define the functor of families in terms of norms and traces
and prove that the two definitions are equivalent when $d!$ is invertible. 
The norm definition
is essentially the one given by D. Rydh in \cite{R1} - in particular the
corresponding functor is represented by the space of divided powers -
while the trace definition is a generalization of the one given in 
\cite{An}. We also 
obtain a formula for the tangent space to a point in the symmetric power,
which appears to be new.
  In Section 4 we prove that divided powers of a line bundle give a line 
bundle, define norm functors and use them to formulate
 a third definition for families of zero cycles. We 
 prove that it is equivalent to the definition via  norm
maps. Finally, in Section 5 we interpret 
in terms of norm functors such standard constructions as 
Hilbert-Chow morphism, sums and direct images of 
cycles, and Chow forms. Quite naturally, our descriptions are closely
related to those of \cite{R2}.

\bigskip
\noindent
\textbf{Acknowledgements.} This work was supported by the Sloan Research 
Fellowship. The author is also grateful to V. Drinfeld from whom he 
learned Grothendieck's approach to zero cycles and to D. Rydh for his useful
remarks.

\section{Preliminaries}

\subsection{Polynomial laws and divided powers.}

We recall some definitions from \cite{Ro}, cf. also 
\cite{La}. In this subsection all rings and algebras will be assumed
commutative and with unity, although the theory 
can be developed in greater generality, cf. 
\textit{loc. cit.}
Let  $M, N$ be two modules over a ring $A$. 
Denote by $\mathcal{F}_M$ the functor 
$$
\mathcal{F}_M: A\textrm{-alg} \to Sets, \qquad A' \to A' \otimes_A M
$$
where $A\textrm{-alg}$ is the category of (commutative) $A$-algebras. 

\bigskip
\noindent
\textbf{Definition} 
A \textit{polynomial law} from $M$ to $N$ 
is a natural transformation $F: \mathcal{F}_M \to \mathcal{F}_N$, 
i.e. for every $A$-algebra $A'$ it defines 
a map $F_{A'}: A' \otimes_A M \to A' \otimes_A N$ and for any morphism 
$A' \to A''$ of $A$-algebras the natural agreement condition is 
satisfied. The polynomial law $F$ is \textit{homogeneous of degree $d$} 
is $F_{A'}(ax) = a^d F_{A'}(x)$ for any $a \in A'$ and 
$x \in A' \otimes_A M$. 
If $B$ and $C$ are $A$-algebras, then $F: \mathcal{F}_B
\to \mathcal{F}_C$ is \textit{multiplicative} if 
$F_{A'}(1) = 1$ and $F_{A'} (xy) = F_{A'}(x) F_{A'}(y)$ for $x, y 
\in B \otimes_A A'$. 

\bigskip
\noindent
Denote by $Pol^d(M, N)$ the set of homogeneous polynomial laws of
degree $d$. 
By Thm IV.1 on p. 266 in \cite{Ro}, the functor $Pol^d(M, ?)$ is 
representable: there exists an $A$-module $\Gamma^d_A(M)$, 
called \textit{module of degree $d$ divided powers}, 
and an isomorphism of functors in $N$: 
\begin{equation}
\label{roby}
Pol^d(M, N)  \simeq  Hom_A(\Gamma^d_A(M), N) 
\end{equation}
Moreover, if $B, C$ are $A$-algebras then each 
$\Gamma^d_A(B)$ is also an $A$-algebra and multiplicative 
laws in $Pol^d(B, C)$ correspond precisely to $A$-algebra 
morphisms $\Gamma^d_A(B) \to C$, cf. Theorem 7.11 in 
\cite{La} or Proposition 2.5.1 in \cite{Fe}. 

\bigskip
\noindent
Explicitly, the direct sum $\Gamma_A(M)= \bigoplus_{d \geq 0} \Gamma^d_A(M)$
may be defined as a unital graded commutative  $A$-algebra with product
$\times$, degree $d$ generators $\gamma^d(x)$, $x \in M, d \geq 0$  and relations
$$
\gamma^0(x) = 1; \qquad \gamma^d(xa) = \gamma^d(x) a^d; 
\qquad \gamma^d(x) \times \gamma^e(x) = \binom{d+e}{e} \gamma^{d+e}  (x);
$$
$$
\gamma^d(x + y) = \sum_{d_1 + d_2 = d} \gamma^{d_1}(x) \times \gamma^{d_2}(y)
$$
In particular, $\Gamma^0_A(M)\simeq A$ and $\Gamma^1_A(M)\simeq M$
with $\gamma^1(x)$ given by $x$. 
We briefly summarize the properties of this construction
\begin{enumerate}
\item $\Gamma_A(\cdot)$ is a covariant functor from the 
category of $A$-modules to the category of graded $A$-algebras which 
commutes with base change $A \to A'$. 

\item If $B$ is an $A$-algebra, then the $A$-algebra $\Gamma^d_A(B)$
satisfies $\gamma^d(xy) = \gamma^d(x) \gamma^d(y)$ for any $x, y \in B$. 
Below we will also use a formula for arbitrary products which can be
found in 2.4.2 of \cite{Fe}.

\item The map $\gamma^d: M \to \Gamma^d_A(M)$ is a homogeneous polynomial law
of degree $d$. The isomorphism of \eqref{roby} is obtained by composing
an $A$-module homomorphism $\Psi_n:\Gamma^d_A(M)\to N$ with $\gamma^d$
to obtain a polynomial law $n: M\to N$:
$$
n = \Psi_n \circ \gamma^d
$$

\item When $M$ is flat over $A$ or $d!$ is invertible in $A$, 
$\Gamma_A(M)$ is isomorphic to 
the algebra of symmetric tensors $TS_A(M)$, i.e. the subalgebra 
$\oplus_{d \geq 0} \big[T^d_A(B)\big]^{\Sigma_d}$ in the 
tensor algebra $T_A(B)$ equipped with the commutative shuffle product.  
In the second case we can further identify both algebras with the symmetric 
algebra $S_A(M)$ (i.e. the \textit{quotient} of the tensor algebra 
$T_A(M)$ by the obvious relations). 

\item If $F \in Pol^d(M, N)$ and we evaluate $F$ at the $A$-algebra
$A' = A[t_1, \ldots t_k]$ then $F(t_1 m_1 + \ldots t_k m_d)
\in N[t_1, \ldots, t_k]$ is a sum of degree $d$ monomials in 
$t_1, \ldots, t_k$ and the coefficient of $t_1^{\alpha_1} \ldots 
t_k^{\alpha_k}$ is the value of the corresponding $A$-module  
homomorphism $\Psi_{F}: \Gamma^d_A(M)\to N$ at 
$\gamma^{\alpha_1}(x_1) \times \ldots \times \gamma^{\alpha_k} (x_k)$. 
This explains the term ``degree $d$ homogeneous polynomial law".

\end{enumerate}

\subsection{Norms and traces for finite flat morphisms.}

Let $\pi: Y \to S$ be a finite flat morphism of schemes or algebraic
spaces and assume that 
$\pi_* \mathcal{O}_Y$ is locally free of constant rank $d$. We 
have a natural morphism of $\mathcal{O}_S$-modules 
$$
\pi_* \mathcal{O}_Y \to End_{\mathcal{O}_S} (\pi_* \mathcal{O}_Y)
$$
and taking the composition with trace and determinant 
we obtain two maps
$$
\theta: \pi_* \mathcal{O}_Y \to \mathcal{O}_S; 
\qquad n: \pi_* \mathcal{O}_Y\to \mathcal{O}_S
$$
It is easy to see that $\theta$ is a morphism of $\mathcal{O}_S$-modules
and that $n$ extends to a homogeneous polynomial law of degree $d$ and
therefore defines a section $\sigma: S \to Spec(\Gamma^d_{\mathcal{O}_S}
(\mathcal{O}_Y))=:\Gamma^d(Y/S)$.
For any line bundle $L$ on $Y$ we also define its \textit{norm}
$$
N(L) = Hom_{\mathcal{O}_S}(\Lambda^d(\pi_* \mathcal{O}_Y), 
\Lambda^d(\pi_* L)).
$$
On the other hand, if $S = Spec(A)$, $Y= Spec(B)$ are affine
and $L$ is given by an invertible $B$-module $M$ then $\Gamma^d_A(M)$
is naturally a $\Gamma^d_A(B)$-module, cf. \cite{Ro}, \cite{La}. 
Using this construction locally on $S$ in the general case, we
obtain a quasi-coherent sheaf on $\Gamma^d(Y/S)$ which we 
denote by $\Gamma^d(L)$. One can show that $\Gamma^d(L)$ is 
an invertible $\mathcal{O}$-module whenever $L$ is, see Section 4.1 below.
We have the following important result
\begin{lemma}
In the notation introduced above, $N(L)\simeq \sigma^* \Gamma^d(L)$.
We also have canonical isomorphisms
$$
N(L \otimes_{\mathcal{O}_Y} F) \simeq N(L) \otimes_{\mathcal{O}_S} N(F); 
\qquad N(\pi^*(L)) \simeq L^{\otimes d}
$$
where $F$ is an invertible $\mathcal{O}_Y$-module.
\end{lemma}
\textit{Proof} The first two assertions are proved
in  Proposition 3.3 in \cite{Fe} (in fact, $F$ can be a coherent sheaf on $Y$ 
if the norm is understood as in \text{loc. cit}), and the third follows from 
the above definition of $N(L)$ and the projection formula. $\square$

\section{Families of zero cycles via norm and trace maps}

\subsection{Relation between traces and norms.}

First assume that $S = Spec(A)$, $X = Spec(B)$ and $d!$ is 
invertible in $A$
and fix a homogeneous polynomial law $n: B\to A$ of degree $d$. 
Then $\Gamma^d_A(B)$ is
isomorphic to the algebra of symmetric tensors
$TS^d_A(B)$, cf. \cite{La}. On one hand,
$A$-algebra homomorphisms $\Gamma^d_A(B) \to A$ correspond to 
homogeneous multiplicative polynomial laws $B\to A$ of degree $d$.
On the other hand,   homomorphisms $TS^d_A(B)$ are  described by 
certain ``trace morphisms" $\theta: B \to A$ of $A$-modules, 
cf.  \cite{An} and \cite{BR} (in the latter paper they are called
Frobenius $n$-homomorphisms).
We briefly outline the  relation between polynomial laws and 
trace morphisms. It does not seem to appear in the literature as
explicitly as below, although many ingredients can be found in 
\cite{An}, \cite{BR} and in Iversen's formalism of linear determinants,
cf. \cite{Iv}. 

The polynomial law $n$ gives, in particular, a map
$n_{A[t]}: B[t] \to A[t]$. Imitating the relationship between 
determinant and trace of a linear operator (compare 
also with Section 2.2) we define $\theta: B \to A$ as the coefficient of $t$ in 
$n_{A[t]}(1 + b t)$. More generally, for $k \geq 1$ define a map
$$
\Theta_k: B^{\times k} \to A
$$
by sending $(b_1, \ldots, b_k)$ to
the coefficient of $t_1 \ldots t_k$ in 
$n_{A[t_1, \ldots, t_k]}(1 + t_1 b_1 + \ldots + t_k b_k)$. 
In terms of the morphism of $A$-algebras $\Psi_n: \Gamma^d_A(B)
\to A$ corresponding to $n$, by property (5) in Section 2.1 we
have
$$
\Theta_k(b_1, \ldots, b_k) = 
\Psi_n(\gamma^{d-k} (1) 
\times  \gamma^1(b_1) \times \ldots \times \gamma^1(b_k))
$$
\begin{lemma}
\label{trace}
For any degree $d$ polynomial law $n$ the
maps $\Theta_k$, $k \geq 1$ satisfy the following properties
\begin{enumerate}
\item $\Theta_k = 0$ for $k > d$;
\item $\Theta_1 = \theta$ and $\Theta_d(x, \ldots, x) = d! n(x)$;
\item $\Theta_k$ is  symmetric in 
its arguments and $A$-linear in each of them, i.e.
descends to an $A$-module morphism from the $k$-th
symmetric power $S^k_A(B) \to A$;
\item if, in addition, $n$ is multiplicative then
the following formula holds for all $k \geq 1$
$$
\Theta_{k+1}(b_1, \ldots, b_{k+1}) :=
$$
$$
\theta(b_1) \Theta_k(b_2, \ldots, b_{k+1}) - 
\Theta_k(b_1 b_2, b_3, \ldots, b_{k+1}) -
\Theta_k(b_2, b_1 b_3, \ldots, b_{k+1}) - \ldots -
\Theta_k(b_2, b_3,  \ldots, b_1 b_{k+1}) 
$$
\end{enumerate}
\end{lemma}
\textit{Proof.} The properties (1), (2) immediately follow
from the definitions and the identity $[\gamma^1(x)]^{\times d} = 
d! \gamma^d(x)$. 
In (3), symmetry also follows immediately from the definitions.
Part (4) follows from multiplicativity of $\Psi_n$ and the product
formula, cf. 2.4.2 in \cite{Fe}:
$$
[\gamma^{d-1} (1) \times \gamma^1(b_1)] 
[\gamma^{d-k} (1) \times \gamma^1(b_2) \times \ldots
\times \gamma^1(b_{k+1})] =
$$
$$
\gamma^{d-k-1}(1) \times \gamma^1(b_1)\times \ldots \times \gamma^1(b_{k+1})
+ \sum_{i =2}^{k+1} \gamma^{d-k}(1) \times \gamma^1(b_2) \times
\ldots \times \gamma^1(b_1 b_i) \times \ldots \times \gamma^1(b_{k+1}).
$$
Finally, $A$-multilinearity in (3) 
follows from the linearity of $\gamma^1$ and the linearity of $\Psi_n$. 
$\square$

\bigskip
\noindent
\textbf{Definition.} Let $B$ be an $A$-algebra. A morphism of $A$-modules
$\theta: B\to A$ is a \textit{degree $d$ trace} if
$$
\theta(1) = d; \qquad \Theta_{d+1} \equiv 0. 
$$
where $\Theta_k$ are constructed from $\theta = : \Theta_1$ using
the formula (4) in Lemma \ref{trace}.

\bigskip
\noindent
\textbf{Remark.} Since $\theta: B \to A$ admits an 
obvious $A'$-linear extension to $\theta_{A'}: A' \otimes_A B
\to A'$ for any $A'$-algebra $A$, 
applying part (2) of the lemma, we see that
the polynomial law $n$ can be recovered from $\theta$ completely. 
Conversely, given a degree $d$  trace $\theta: B\to A$, 
one can check that $n(x) = \frac{1}{d!} \Theta_d(x, \ldots, x)$
defines a polynomial law $n$. In fact, $n(x)$ is homogeneous
of degree $d$ by part (3) of Lemma \ref{trace} and 
multiplicative by Theorem 1.5.3 in \cite{An} or Theorem 2.8 in \cite{BR}. 
Since $\theta$
has canonical pullbacks $\theta \otimes 1:B\otimes A' \to A'$
for all $A$-algebras $A'$, so does $n$, i.e. we obtain 
a polynomial law.

Observe  that we can also define $\theta$ as a map ``tangent"
to $n$, i.e. by considering the $A$-algebra $A' = A[\varepsilon]/\varepsilon^2$
and then using the identity
$$
n_{A'} (1 + \varepsilon b) = 1 + \varepsilon \theta(b).
$$
\begin{lemma}
\label{bijections}
The operations $\theta(b) \mapsto n(b) = \frac{1}{d!}\Theta_d(b, \ldots, b)$
and $n(b) \mapsto \theta(b) = \Psi_n(\gamma^{d-1}(1) \times \gamma^1(b))$
define mutually inverse bijections between the set of degree $d$ traces
and the set of degree $d$ norm maps.  
\end{lemma}
\textit{Proof.} We have seen before that traces are sent to norm maps and the
other way around. Let us show that the two constructions are
mutually inverse to each other. 

First assume that $n(x) = \frac{1}{d!} \Theta_d(x, \ldots, x)$
and let us show that the trace constructed from $n$ coincides with 
the original $\theta = \Theta_1$.
If $\Theta_k$ are defined from $\theta$ using formula (4) in
Lemma \ref{trace} one can show that $\Theta_k$ are symmetric
and multilinear (see e.g. Definition 1.3.1 in \cite{An} where
$\Theta_k$ are denoted by $P^k_\theta$). 
The polynomial law
$n(x)$ gives an $A$-algebra homomorphism $\Psi_n: \Gamma^d_A(B) \to A$ 
and 
$$
\Theta_d(b, \ldots, b) = \Psi_n(d! \gamma^d(b)) = \Psi_n (\gamma^1(b)\times
\ldots \gamma^1(b))
$$
Since $\gamma^1(b)$ is $A$-linear in $b$ and $d!$ is invertible in $A$, 
by an easy polarization argument we conclude that 
$\Theta_d(b_1, \ldots, b_d) = \Psi_n(\gamma^1(b_1) \times \ldots 
\times \gamma^1(b_d))$. Using the recursive definition of $\Theta_k$
we get
$$
\Theta_{k+1}(1, b_2, \ldots, b_{k+1}) 
= (d-k) \Theta_k(b_2, \ldots, b_{k+1})
$$
and by descending induction on $k$ we conclude that
$\Theta_k(b_1, \ldots, b_k)= \Psi_n(\gamma^{d-1}(1) 
\times \gamma^1(b_1) \times \ldots \times \gamma^1 (b_k))$.
In particular, 
$$
\theta(b) = \Psi_{\frac{1}{d!} \Theta_d(b, \ldots, b)}
(\gamma^{d-1}(1) \times \gamma^1(b)).
$$
On the other hand, if we start with a polynomial law $n(x)$ and 
set $\theta(b)= \Psi_n(\gamma^{d-1}(1) \times \gamma^1(b))$
then Lemma \ref{trace} tells us that $n(b) = \frac{1}{d!} 
\Theta_d(b, \ldots, b)$, as required. $\square$

\bigskip
\noindent
\textbf{Examples.} 

\bigskip
\noindent
(1) Let $f_i: B\to A$ be $A$-algebra homomorphisms
for $i = 1, \ldots, d$. Then the product $n = f_1 \ldots f_d$ is a degree $d$
homogeneous polynomial law and $\theta= f_1 + \ldots + f_d$ is
the degree $d$ trace corresponding to it, while $\Theta_k$ is given 
by the $k$-th elementary symmetric function in the $f_i$ (up to a scalar). 

\bigskip
\noindent
(2) Let $A = k$ be a field of characteristic
$p$ with $p=0$ or $p > d$ and $Q$ a $k$-point of $Spec(B)$ with $k(Q) = k$, 
corresponding to the evaluation homomorphism $B\to k$, $b \mapsto b(Q)$.
Consider the polynomial law $b \mapsto b(Q)^d$
corresponding to the effective cycle $[dQ] \in Spec(\Gamma^d_k(B))
\simeq Sym^d(Spec(B))$.
We have the following formula for the dual of the tangent space at
$[dQ]$:
$$
T_{[dQ]}^\vee \simeq \mathfrak{m}_Q/\mathfrak{m}^{d+1}_Q
$$
In fact, by assumption on $k$
an element of $T_{[dQ]}$ corresponds to a degree $d$ trace
$$
\theta = \theta' + \varepsilon \theta'': B\to k[\varepsilon]/\varepsilon^2
= k \oplus \varepsilon k
$$
with $\theta'(f) = d \cdot f(Q)$. Since $\theta(1) = d$, $\theta''$ 
vanishes on the subspace of constants $k \subset B$ and therefore 
we can identify it with a linear function of $\mathfrak{m}_Q$. Let
us show that the condition $\Theta_{d+1} (b_1, \ldots, b_{d+1}) = 0$
for all $b_i \in B$, is equivalent to $\theta''(\mathfrak{m}_Q^{d+1}) = 0$. 
In fact, since $\Theta_{d+1}$ is multilinear, we can assume that each 
$b_i$ is either $1$ or in $\mathfrak{m}_Q$. If at least one of the $b_i$ is
$1$, by symmetry we can assume that $b_1 = 1$ and then $\theta(1) = d$
together with the formula (4) in Lemma \ref{trace} immediately give the 
vanishing of $\Theta_{d+1}$. If all the arguments $b_1, \ldots, b_{d+1}$
are in $\mathfrak{m}_Q$, then it is easy to show by induction
using the same formula, that
$\Theta_{l+1} (b_1, \ldots, b_{l+1}) = (-1)^l l! \;\varepsilon\theta''(b_1 \ldots b_{l+1})
$ with $l \geq 0$ and the usual convention $0!=1$. 
In particular, $\Theta_{d+1} (b_1,\ldots, b_{d+1})= 
(-1)^d d!\; \varepsilon \theta''(b_1 \ldots b_{d+1})$ and 
using our assumption on $k$ again, we see that 
$\Theta_{d+1} \equiv 0$ if and only if $\theta''$ descends
to a linear function on $\mathfrak{m}_Q/\mathfrak{m}_Q^{d+1}$, 
which proves the assertion. 

\subsection{Functors of zero cycles.}

\textbf{Definitions.}

\medskip
\noindent
(1) Let $\pi: X \to S$ be a separated morphism of algebraic spaces, cf. 
\cite{Kn}.
Let $Chow^{n}_{\pi, d}$ be a functor on the category of 
algebraic spaces over $S$, sending $T \to S$ to a set 
$Chow^n_{\pi, d}(T)$ of equivalence classes of pairs
$(Y, n)$, where $Y\hookrightarrow X_T$ is a closed algebraic subspace
which is integral over $T$ (i.e. affine and such that locally over
$T$ every regular function on $Y$ satisfies a monic polynomial with 
coefficients in $\mathcal{O}_T$), 
and $n: (\pi_T)_*\mathcal{O}_Y\to \mathcal{O}_T$
is a multiplicative polynomial law of degree $d$. Two pairs $(Y_1, n_1)$, 
$(Y_2, n_2)$ are called equivalent if there is a third pair $(Y, n)$,
such that $Y$ is an algebraic subspace in $Y_1 \cap Y_2$ which is
integral over $T$, and $n_i$ is equal to the composition of the
natural morphism $(\pi_T)_* \mathcal{O}_{Y_i} \to (\pi_T)_*\mathcal{O}_{Y}$ 
with $n$, for $i = 1, 2$. The inverse image of $(Y, n)$ with respect to 
an $S$-morphism $\phi: T' \to T$ is given by $(Y', n')$ where 
$Y' = Y\times_T T'$ and $n'$ is described on elements of an 
affine covering $T' = \cup Spec(A_i)$ by restricting $n$ to those 
$\mathcal{O}_T$-algebras which factor through $A_i$.

\medskip
\noindent
(2) Let $\pi: X \to S$ be as before and assume that $d!$ defines
an invertible regular function on $S$. Let $Chow^\theta_{\pi, d}$
be a functor on the category of algebraic spaces over $S$, given by 
equivalence classes of pairs $(Y, \theta)$, where $Y \hookrightarrow X_T$
is a closed algebraic subspace of $X_T$ which is integral over $T$ and
$\theta: (\pi_T)_* \mathcal{O}_Y\to \mathcal{O}_T$ is a degree $d$ trace.
Equivalence of such pairs is defined in a similar way. Note that 
for a pullback $(Y', \theta')$ with respect to $\phi: T' \to T$ 
we can define $\theta'$ 
simply as $\theta \otimes_{\mathcal{O}_T} \mathcal{O}_{T'}$.

\begin{prop}
Assume that $d!$ is an invertible regular function on $S$. 
There exists an isomorphism of functors $Chow^n_{\pi, d} \simeq
Chow^\theta_{\pi, d}$.
\end{prop}
\textit{Proof.} Follows immediately from Lemma \ref{bijections}
$\square$

\bigskip
\noindent
\textbf{Remarks.} 

\medskip
\noindent
(1)  The definition of $Chow^{n}_{\pi, d}$ is a 
restatement of Definition 3.1.1 in \cite{R1}. Therefore, 
this functor is represented by the space of divided powers $\Gamma^d(X/S)$. 
The definition of $Chow^{\theta}_{\pi, d}$ is a version of 
Definition on page 7 of \cite{An} applied to zero cycles, but stated
here in greater generality.

\bigskip
\noindent
(2) Let $n: B \to A$ be a degree $d$ norm map. 
Then, following the Definition 2.1.5 in \cite{R1}
we define the characteristic polynomial of $b \in B$ by the formula
$$
\chi_{n, b}(t):= n_{A[t]} (b-t) = \sum_{k=0}^d (-1)^k \Psi_n
(\gamma^k(1) \times \gamma^{d-k}(b)) t^k \in A[t]
$$
In the notation of Lemma \ref{trace} we have
$\chi_{n, b}(t) = \sum_{k=0}^d (-1)^k \frac{\Theta_{d-k}(b, 
\ldots, b)}{(d-k)!}  t^k$. Now let $J_n \subset B$ be the ideal generated 
by $\chi_{n, b}(b)$ for all $b \in B$, called
the \textit{Cayley-Hamilton ideal} of $n$. Then by Proposition 2.1.6
in \cite{R1} the norm map $n$ factors through the quotient $B/J_{n}$. 
Similarly, for a degree $d$ trace $\theta: B \to A$
the Section 1.6 of \cite{An} defines an ideal $J_\theta \subset B$
such that $\theta$ factors through $B/J_{\theta}$. 

We observe here that $J_n = J_\theta$ if $n$ and $\theta$ are 
related by the bijection of Lemma \ref{trace}. In fact, by 
Definition 1.6.2.2 of \textit{loc. cit.}
$J_\theta$ is generated by values of the polarized version of $\chi_{n, b}$
only $\Theta_{d-k}$ are defined based in $\theta$, as in part (4)
of Lemma \ref{trace} while in the case of $J_n$ they are defined through
$n$. Since $d!$ is assumed invertible in $A$ the values of the polarized
version generates the same ideal as values of $\chi_{n, b}$ itself.

\bigskip
\noindent
(3) In \cite{An} traces were defined 
on the completion $\widehat{X}_T$ at
a closed subset $Z \subset X_T$ which is proper and of pure 
relative dimension zero over $T$. But by Corollary 1.6.3
in \textit{loc. cit.},  a degree $d$ 
trace $\theta:(\pi_T)_* \mathcal{O}_{\widehat{X}_T} \to \mathcal{O}_T$ 
descends to 
the subscheme $Y$ given by $\mathcal{O}_{Y'}/\mathcal{J}_\theta$ 
which is integral over $T$. Therefore we can restrict to integral subschemes
in the definition.

\section{Families of zero cycles via norm functors}

\subsection{Divided powers of line bundles.}

Let $\pi: X \to S$ be an affine morphism of algebraic spaces. 
For any quasi-coherent sheaf $L$ on $X$, 
the sheaf $\Gamma^d_{\mathcal{O}_S} (\pi_* L)$ is a module over the 
$\mathcal{O}_S$-algebra $\Gamma^d_{\mathcal{O}_S}  (\pi_*\mathcal{O}_X)$, 
cf. \cite{La}.
This gives a quasi-coherent sheaf $\Gamma^d(L)$ on $\Gamma^d(X/S)$. 

We now recall a construction from Section 3 of \cite{Fe} presenting it
here in a sheafified version. Let $L'$ and $L''$ be two quasi-coherent
sheaves on $X$. There is a unique functorial morphism of $\mathcal{O}_S$-modules
$$
\Gamma^d_{\mathcal{O}_S}(\pi_* L') \otimes_{\mathcal{O}_S} 
\Gamma^d_{\mathcal{O}_S}(\pi_* L'')
\to \Gamma^d_{\mathcal{O}_S}(\pi_* L' \otimes_{\mathcal{O}_S} \pi_* L'')
$$
which sends $\gamma^d(x) \otimes \gamma^d(y)$ to $\gamma^d(x \otimes y)$, 
cf. \cite{Ro}, \cite{La}, \cite{Fe}. The image on a general element is
given by formula 2.4.2 in \cite{Fe}. The composition of this map 
with $\Gamma^d_{\mathcal{O}_S} (\pi_* L' \otimes_{\mathcal{O}_S} \pi_* L'')
\to \Gamma^d_{\mathcal{O}_S}(\pi_* (L' \otimes_{\mathcal{O}_X} L''))$
descends to a morphism of $\Gamma^d_{\mathcal{O}_S} 
(\pi_* \mathcal{O}_X)$-modules
\begin{equation}
\label{tensor}
\Gamma^d_{\mathcal{O}_S}(\pi_* L')
\otimes_{\Gamma^d_{\mathcal{O}_S}(\pi_* \mathcal{O}_X)}
\Gamma^d_{\mathcal{O}_S}(\pi_* L'')
\to \Gamma^d_{\mathcal{O}_S}(\pi_*(L' \otimes_{\mathcal{O}_X} L''))
\end{equation} 
The following result does not seem to appear in the literature:
\begin{lemma}
The morphism \eqref{tensor} is an isomorphism if at least one of the 
sheaves $L'$, $L''$ is an invertible $\mathcal{O}_X$-module. In particular
$\Gamma^d(L)$ is an invertible module on $\Gamma^d(X/S)$ if $L$ is
an invertible module on $X$. The map $L \mapsto \Gamma^d(L)$ 
extends to a functor $N: PIC(X) \to PIC(\Gamma^d(X/S))$ between the
categories $PIC$ of invertible modules (and morphisms as $\mathcal{O}$-modules).
The functor $N$ is equipped with isomorphisms $\Gamma^d(L')
\otimes_{\mathcal{O}_{\Gamma^d(X/S)}} \Gamma^d(L'')
\simeq \Gamma^d(L' \otimes_{\mathcal{O}_X} L'')$ which agree
with commutativity and associativity isomorphisms for tensor product
of line bundles. If $\pi^d: \Gamma^d(X/S) \to S$ is the
canonical morphism then the induced map
$$
\pi_* \mathcal{H}om_{\mathcal{O}_X} (L', L'')
\to \pi^d_* \mathcal{H}om_{\mathcal{O}_{\Gamma^d(X/S)}}
(\Gamma^d(L'), \Gamma^d(L''))
$$
extends canonically to a homogeneous polynomial law of degree $d$. 
\end{lemma}
\textit{Proof.} We here prove the first two assertions, 
since the statements about the $PIC$ functor
follow from a routine check based on the definitions involved.

To prove the isomorphism we can assume 
that $X = Spec(B)$ and 
$S = Spec(A)$ are affine and $L'$ is given by 
a projective $B$-module $P$ of rank 1. Observe
that for any finite subset of points $x_1, \ldots, x_l$ there is an 
affine open subset $U \subset X$ containing these points,
such that $L|_U$ is trivial. In fact,
we can assume that no  $x_i$ is in the closure of another $x_j$
(otherwise we can erase $x_j$ from the list, since trivialization of
$L$ in a neighborhood of $x_i$ will also give a trivialization for $x_j$). 
Since $X$ is affine, we can choose sections $l_1, \ldots, l_d$ in $H^0(X, L)$,
generating the stalks $L_{x_1}, \ldots, L_{x_d}$, respectively. Also, we
can choose functions $f_1, \ldots, f_d$ in $H^0(X, \mathcal{O}_X)$ such 
that each (image of) $f_i$ generates the stalk  $\mathcal{O}_{x_i}$ and
vanishes at $x_j$ for $i \neq j$. Then the section $l = l_1 f_1 + \ldots 
+ l_d f_d$ generates the stalks of $L$ at $x_1, \ldots, x_d$ and
hence defines a trivialization of $L$ in an affine neighborhood $U$
of $x_1, \ldots, x_d$. 

Choose and fix a point $\beta \in \Gamma^d(X/S)$. It suffices to prove that
\eqref{tensor} is an isomorphism in a neighborhood of $\alpha$. 
By the results of Section 2 in \cite{R1} one can find finitely many points
$x_1, \ldots, x_l$ in $X$, with $l \leq d$ and a closed affine
subscheme $Y \subset X$ supported on the union of $x_i$, such 
that $\beta$ is a point in the closed subscheme $\Gamma^d(Y/S)$. 
Choosing an open affine neighborhood $U$ of $x_1, \ldots, x_l$ and
a trivialization $L'|_U$ as above, we obtain a open affine 
neighborhood $\Gamma^d(U/S) \subset \Gamma^d(X/S)$ of $\beta$, and the
trivialization of $L'$ on $U$ induces an isomorphism of $\Gamma^d(L')|_{\Gamma^d(U/S)}$ with the structure sheaf. Thus, on 
$\Gamma^d(U/S)$ the map \eqref{tensor} is an isomorphism and 
$\Gamma^d(L')$ is locally free of rank 1. Note that, once we know 
the isomorphism, the fact that $\Gamma^d(L')$ is invertible 
may also be proved by choosing $L''$ to be the dual of $L'$.
$\square$

\bigskip
\noindent
\textbf{Remarks.} 

\bigskip
\noindent
(1) By equation (2.4.3.1) of \cite{Fe} for any invertible 
$\mathcal{O}_S$-module $P$ and any $\mathcal{O}_X$-module $F$ one has
$$
\Gamma^d(\pi^* P \otimes_{\mathcal{O}_X} F) 
\simeq (\pi^d)^*(P^{\otimes d}) \otimes_{\mathcal{O}_{\Gamma^d(X/S)}}
\Gamma^d(F)
$$

\noindent
(2) Suppose that $S$ can be covered by affine open subsets
$V_i$ such that $L$ is trivial on $\pi^{-1}(V_i)$ (e.g. 
that we are in the situation of Section 2.2). Then 
$\Gamma^d(L)$ is trivial on the open subset $(\pi^d)^{-1}(V_i)$. 
If $\phi_{ij}$ are the transition functions for $L$ then
their norms $\gamma^d(\phi_{ij})$
are the transition functions of $\Gamma^d(L)$. This construction 
was originally given for the setting of Section 2 by Grothendieck in
\cite{EGA2}, 6.5. 

\subsection{Norm functors}

\textbf{Definitions.} 

\medskip
\noindent
(1) Let $\pi: Y \to S$ be a separable morphism of algebraic spaces.
Denote by $\mathcal{P}IC(Y/S)$ the category with objects given by $(S', L)$
where $S' \to S$ is an algebraic space over $S$
and $L$ is a line bundle on $Y_{S'} = Y\times_S S'$. A morphism 
$(\xi, \rho): (S_1, L_1) \to (S_2, L_2)$ in the category 
$\mathcal{P}IC(Y/S)$ is given by a morphism 
$\xi: S_1 \to S_2$ of algebraic spaces over $S$, plus a morphism  
$\rho: L_1 \to \xi^*(L_2)$ of coherent sheaves on $Y_{S_1}$. 
There is an obvious forgetful functor $p_Y: \mathcal{P}IC(Y/S) \to \textrm{Sp}/S$
to the category of algebraic spaces over $S$, given by 
$(T, L)\mapsto T$ and $(\xi, \rho)\mapsto \xi$. 
When $Y= S$ and $\pi$ is the identity morphism we write $\mathcal{P}IC(S)$
instead of $\mathcal{P}IC(S/S)$. There is a natural functor $\mathcal{P}IC(S)
\to \mathcal{P}IC(Y/S)$ given by pullback of $L$ from $S'$ to $Y_{S'}$. We
denote this functor simply by $\pi^*$. 

\medskip
\noindent
(2) Let $\pi: Y \to S$ be as above. A \textit{norm functor of degree $d$
over $\pi$} is a triple $\mathcal{N}= (N, \mu, \epsilon)$ where 
$N$ is a functor 
$\mathcal{P}IC(Y/S) \to \mathcal{P}IC(S)$ such that $p_S \circ N = 
p_Y$. In other words, a pair
$(S', L)\in Ob(\mathcal{P}IC(Y/S))$ is sent to a pair $(S', M) \in 
Ob(\mathcal{P}IC(S))$ and sometimes we will abuse notation by dropping $S'$
and writing $M= N(L)$. Further, 
 for any pair $(S', L_1)$, $(S', L_2)$ of objects in $\mathcal{P}IC(Y/S)$
with the same $S'$ we require an isomorphism 
$$
\mu_{S', L_1, L_2}: N(L_1)\otimes_{\mathcal{O}_{S'}} N(L_2)\simeq 
N(L_1 \otimes_{\mathcal{O}_{Y_{S'}}} L_2)
$$
such that the system of isomorphisms $\mu= \mu_{\{S', \cdot, \cdot\}}$ 
agrees with the the base cange and the standard
symmetry and associativity isomorphisms for line bundles
on $Y_{S'}$ and $S'$, respectively
(see e.g. the last two diagrams on p. 36 of \cite{Du}). 
Finally $\epsilon$ is an isomorphism of functors $\mathcal{P}IC(S) \to
\mathcal{P}IC(S)$:
$$
\epsilon: N\circ \pi^* \simeq (\;\cdot\;)^{\otimes d}
$$
such that $\mu_{\{S', \cdot, \cdot\}}\circ (N \pi^* \otimes N\pi^*)$ is given by 
the canonical isomorphism $L_1^{\otimes d} \otimes_{\mathcal{O}_{S'}} L_2^{\otimes d} \simeq
(L_1 \otimes_{\mathcal{O}_{S'}} L_2)^{\otimes d}$.

\bigskip
\noindent
(3) Let $\pi: X \to S$ be a separated morphism of algebraic spaces. 
Let $Chow^N_{\pi, d}$ be the functor 
from the category $Sp/S$ of algebraic spaces over $S$ to sets, 
sending
$T \to S$ to equivalence classes of data $(Y, \mathcal{N})$
where $Y \hookrightarrow X_T$ is a closed algebraic subspace, 
integral over $T$, and $\mathcal{N}$
 is a degree $d$ norm functor over $(\pi_T)|_Y$.  
Two pairs $(Y_1, \mathcal{N}_1)$
and $(Y_2, \mathcal{N}_2)$ are called equivalent if there is a
third subspace $Y\subset Y_1 \cap Y_2$ and a degree $d$ norm functor 
$\mathcal{N}$ over $(\pi_T)|_Y$  together with isomorphisms
between $N_i$ and the composition of $N$ with the restriction 
from $Y_i$ to $Y$, which are further required to agree with $\epsilon_i$
and $\mu_i$ in the obvious sense.

\bigskip
\noindent
\textbf{Remark.} Since by definition a norm functor is
local on $S$, we obtain a map
\begin{equation}
\label{polyn}
(\pi_{S'})_* \mathcal{H}om_{\mathcal{O}_{Y_{S'}}}(L_1, L_2) \to \mathcal{H}om_{\mathcal{O}_{S'}}(N(L_1), N(L_2))  
\end{equation}
This map is not $\mathcal{O}_{S'}$ linear; it is rather a polynomial 
law of degree $d$. To show this, it suffices to assume 
$L_1 = \mathcal{O}_{Y_{S'}}$. In fact, since by definition $N$
preserves tensor products and sends the trivial bundle to the trivial bundle
we have $N(L^\vee) = N(L)^\vee$. The left hand side can be 
rewritten as 
$(\pi_{S'})_*\mathcal{H}om_{\mathcal{O}_{Y_{S'}}}(\mathcal{O}_{Y_{S'}}, L_1^\vee
\otimes L_2)$ while the right hand side becomes 
$$
\mathcal{H}om_{\mathcal{O}_{S'}}(\mathcal{O}_{S'}, N(L_1)^\vee \otimes N(L_2))
\simeq   \mathcal{H}om_{\mathcal{O}_{S'}}(\mathcal{O}_{S'}, 
N(L_1^\vee \otimes L_2))
$$
A local section $f$ of $\mathcal{O}_{S'}$ acts 
on $(\pi_{S'})_*\mathcal{H}om_{\mathcal{O}_{Y_{S'}}}(\mathcal{O}_{Y_{S'}}, L_1^\vee
\otimes L_2)$
by  composition with the ``multiplication by $f$" endomorphism 
of $\mathcal{O}_{Y_{S'}} \simeq 
\pi_{S'}^* \mathcal{O}_{S'}$. By definition of $\epsilon$, the norm
functor sends it  to multiplication by $f^d$. 

\begin{lemma}
Let $\pi: Y \to S$ be an integral morphism of algebraic spaces 
with universally topologically finite fibers and let $L$
be a line bundle on $Y$. Any point $s \in S$ has an etale
neighborhood $U \subset S$ such that the restriction on $L$
on $\pi^{-1}(U)$ is trivial.
\end{lemma}
\textit{Proof.} It suffices to assume that $S$, and hence also $Y$, are affine.
Since the fiber $\pi^{-1}(s)$ is finite, repeating the argument
in Lemma 5 we can find a section $l$ of $L$ on $Y$ which generates the stalks
of $L$ at each of the points in $\pi^{-1}(s)$. The subset $W \subset Y$ 
of points where $l$ \textit{fails} to generate the stalk of $L$ is closed
in $Y$ and disoint from the fiber $\pi^{-1}(s)$. Its image $\pi(W)$ is 
closed in $Y$ since $\pi$ is integral, and does not contain $s$. Hence
$s$ admits an affine Zariski neighborhood $U \subset (Y \setminus \pi(W))$
such that on $\pi^{-1}(U)$ the line bundle $L$ is trivialized by the
section $l$. $\square$

\begin{lemma}
In $\pi: Y\to S$ is integral, 
then any norm functor has no nontrivial automorphisms.
\end{lemma}
\textit{Proof.} A functor automorphism is given 
by a family of isomorphisms $\phi_{(T, L)}: N(L) \to N(L)$
for all objects $(T, L)$ of $\mathcal{P}IC(Y)$. If $L$ is pulled
back from $T$ this automorphism has to be identity since it has
to respect $\epsilon$. By the previous lemma, we can find an etale
open cover $\{U_i\}$ of $T$, such that $L$ is trivial over each $U_i$. 
Then the restriction of $\phi_{(T, L)}$ to each $U_i$ is the identity 
due to the agreement with $\epsilon$ hence  $\phi_{(T, L)}$ is itself identity.
$\square$

\bigskip
\noindent
\textbf{Remark.} For a general $\pi$ the previous result fails. 
One possible example is the situation when
$Y$ and $S$ are over a field $k$, $Y = Y_0 \times_{Spec{k}} S$
and there exists a non-trivial group homomorphism 
$Pic(Y_0) \to \mathcal{O}_S^*$ (where $Pic$ is the group of isomorphism
classes of line bundles on $Y_0$).

\begin{prop}
The functor $Chow^N_{\pi, d}$ is isomorphic to $Chow^n_{\pi, d}$
and is therefore represented by the space of 
divided powers $\Gamma^d(X/S)$. 
\end{prop}
\textbf{Proof.} An $S$-morphism $T \to \Gamma^d(X/S)$ induces a 
norm functor by Lemma 5 by using the pullback of line bundles $\Gamma^d(L)$
to $T$.
Conversely, taking $L_1 = L_2 = \mathcal{O}_Y$ in \eqref{polyn}
we obtain a norm map $(\pi_T)_* \mathcal{O}_{Y_T} \to \mathcal{O}_T$.

It is obvious that $Chow^n_{\pi, d} \to Chow^N_{\pi, d} 
\to Chow^n_{\pi, d}$ is identity since we are essentially expanding the
data involved in the definition of $Chow^n_{\pi, d}$ and then 
forgetting the extra data constructed. 

In the opposite direction, suppose we have a closed subspace $Y \subset
X_T$ and a norm functor $\mathcal{N} = (N, \mu,\epsilon)$ over 
$(\pi_T)|_Y$ which we
use to extract the polynomial law $(\pi_T)_* \mathcal{O}_Y \to \mathcal{O}_T$
and thus obtain an $S$-morphism $\sigma: T \to \Gamma^d(Y/S)$. We need 
to construct isomorphisms $N(L) \simeq \sigma^* (\Gamma^d(L))$ for
all line bundles $L$, which commute with pullbacks, agree with multiplicativity 
isomorphisms and give identity on $L^{\otimes d}$ on bundles pulled back
from $T$ to $Y$. In other words, we need to prove an isomorphism 
$$
\mathcal{O}_S
\otimes_{\Gamma^d_{\mathcal{O}_S} ((\pi_T)_* \mathcal{O}_Y)}
\Gamma^d_{\mathcal{O}_S} ((\pi_T)_* L)   \simeq N(L)
$$
It suffices to construct a morphism from the left hand side to the
right hand side and then apply Lemma 5 to find a Zariski open covering
$\{U_i\}$ of $T$ such that $L$ is trivial on the preimage of $U_i$, in which 
case the isomorphism becomes a tautology. 
To that end, observe that \eqref{polyn} gives a polynomial law
$
(\pi_T)_* L \to N(L)
$
hence a morphism of $\mathcal{O}_S$-modules
$$
\mu_L: \Gamma^d_{\mathcal{O}_S}((\pi_T)_* L) \to N(L)
$$
The fact that $\mu_L$ descends to the above tensor product is equivalent to 
the formula
\begin{equation}
\label{action}
\mu_L(f \cdot s) = \mu_{\mathcal{O}_Y}(f) \cdot \mu_L(s)
\end{equation}
where $f$, resp. $s$, is a local section of 
$\Gamma^d_{\mathcal{O}_S} ((\pi_T)_* \mathcal{O}_Y)$, resp.
$\Gamma^d_{\mathcal{O}_S} ((\pi_T)_* L)$, and the module structure of the
left hand side is given e.g. by formula 7.6.1 in \cite{La}. But \eqref{action}
follows from the fact that $N$ is a functor, i.e. maps compositions of 
morphisms to compositions of morphisms, and the fact that after a faithfully 
flat base change the $\mathcal{O}_T$-module $\Gamma^d_{\mathcal{O}_T} ((\pi_T)_* L)$ is generated by $\gamma^d((\pi_T)_* L)$.

Agreement of $\sigma^*(\Gamma^d(L)) \simeq N(L)$ with multiplicativity isomorphisms and normalization on line bundles pulled back from $T$, also 
follows from the functor property of $N$.
$\square$

\subsection{Non-homonegeous norm functors}

One can also give a definition of a non-homogeneous norm functor as
a triple $(N, \mu, \epsilon)$ where $N$ and $\mu$ are as before and
$\epsilon$ is an isomorphism
$$
\epsilon: N(\mathcal{O}_Y) \simeq \mathcal{O}_S
$$
which sends the identity endomorphism of $\mathcal{O}_Y$ to the identity 
endomorphism of $\mathcal{O}_S$. 

Observe that the proof of Lemma 6 still works in this case and hence
non-homogeneous norm functors form a set. As in the homogeneous case, 
any such functor gives a multiplicative polynomial law
$$
\pi_* \mathcal{O}_Y \to \mathcal{O}_S
$$
and hence by Section 2 of \cite{R2} it defines a section
$$
S \to \Gamma^\star (Y/S) = \coprod_{d \geq 0} \Gamma^d(Y/S).
$$
Repeating the argument of the previous subsection one shows that
the functor of zero cycles $Chow^N_{\pi, \star}$
defined via non-honogeneous norm functors is 
isomorphic to the functor of zero cycles defined via non-homogeneous norm 
maps $Chow^n_{\pi, \star}$. Therefore $Chow^N_{\pi, \star}$ is represented
by the space of effective zero cycles $\Gamma^\star(X/S)$. Details are left 
to the motivated reader. 

\section{Standard constructions}

\bigskip
\noindent
\textbf{Hilbert-Chow morphism.}

\bigskip
\noindent
If $\pi: Y \to S$ finite and flat and $\pi_* \mathcal{O}_Y$
is locally free of constant rank $d$ on $S$ the construction 
of section 2.2 gives the norm of a line bundle $L$ on $Y$. 
Lemma 1 shows that the norm of line bundles defines
a norm functor $\mathcal{P}IC(Y/S) \to \mathcal{P}IC(S)$ inducing a morphism of representing spaces:
$$
Hilb^d(X/S) \to Chow^N_{\pi, d}(X/S).
$$
\textbf{Sums of cycles}

\bigskip
\noindent
If $(Y_1, \mathcal{N}_1)$ and $(Y_2, \mathcal{N}_2)$ are two families of zero 
cycles of
degrees $d_1$ and $d_2$, respectively, then their sum has degree $d_1 + d_2$ 
and is given by $(Y_1 \cup Y_2, (\mathcal{N}_1 \circ i_1^*) \otimes 
(\mathcal{N}_2 \circ i_2^*))$ where $(i_l)_*$ is the functor defined by 
restriction of line bundles from $Y = Y_1 \cup Y_2$ to $Y_l$ for $l = 1, 2$.
This induces the sum morphism
$$
\pi_{d_1, d_2}: \Gamma^{d_1}(X/S) \times_S \Gamma^{d_2}(X/S)
\to \Gamma^{d_1 + d_2} (X/S) 
$$

\bigskip
\noindent
\textbf{Universal family.}

\bigskip
\noindent
For $T = \Gamma^d(X/S)$ consider the base change morphism $\pi_T:
X_T \to T$. Set $Y = \Gamma^{d-1}(X/S) \times_S X$ which maps to $T$ 
via the addition morphism $\pi_{d-1, 1}$. 
By \cite{R2}, $Y$ is integral over $T$ and can be identified with a
closed subspace of $\Gamma^d(X/S) \times_S X$ via the 
morphism $(\xi, x)\mapsto (\xi + x, x)$. There does not seem to be a
quick definition of the corresponding universal norm functor 
$N: \mathcal{P}IC(Y/T) \to \mathcal{P}IC(T)$, as is also the case with 
the universal norm map $(\pi_{d-1, 1})_* \mathcal{O}_Y \to \mathcal{O}_T$. 
However, if $\eta: \Gamma^{d-1}(X/S)\times_S X \to X$ is the canonical 
projection, it follows easily that the composition 
$$
\mathcal{P}IC(X) 
\stackrel{\eta^*}\to \mathcal{P}IC(Y/T) 
\stackrel{N}\to \mathcal{P}IC(T) = \mathcal{P}IC
(\Gamma^d(X/S))
$$
is given simply by the functor $L \mapsto \Gamma^d(L)$.

\bigskip
\noindent
\textbf{Direct image of cycles.}

\bigskip
\noindent
Let $\pi': Z \to S$ be another separated morphism of algebraic spaces
and  $f: X \to Z$ is a morphism over $S$. Take
a faimily of zero cycles on $X$ is represented by a pair $(Y, \mathcal{N})$.
By Section 2 of \cite{R1} we can assume that $Y$ is has univerally 
topologically finite fibers over $S$ and hence by the 
appendix to \textit{loc. cit.} $f(Y)$ is a well-defined
algebraic subspace of $Z$ which is integral over $S$.
One can also give a more direct proof of this result using the 
approximation results (Theorem D) in \cite{R4}. 
The direct image cycle is defined by  $(f(Y), \mathcal{N} \circ f^*)$, which 
induces a morphism 
$$
Chow^N_{\pi, d} (X/S) \to Chow^N_{\pi', d}(Z/S)
$$

\bigskip
\noindent
\textbf{Chow forms.}

\bigskip
\noindent
Assume that $X= Proj(\mathcal{A})$ where $\mathcal{A}
= \bigoplus_{l \geq 0} \mathcal{A}_l$ is a graded $\mathcal{O}_S$-algebra
generated over $\mathcal{A}_0 = \mathcal{O}_S$  by its first component
$\mathcal{A}_1$. Then the natural sheaf $\mathcal{O}(1)$ on $X$ is invertible.

Let $(Y, \mathcal{N})$ be a pair representing an element in 
$Chow^N_{\pi, d}(T)$ with $\xi: T \to S$
and denote the inverse image of $\mathcal{O}(1)$ on $Y$ by 
$L$. By assumption a local section of 
$\xi^* \mathcal{A}_l$ on $U \subset T$ gives
a section of $L^{\otimes l}$ on $\pi^{-1}_T (U)$ and hence by the Remark in 
Section 4.2 a section of $N(L^{\otimes l}) 
\simeq N(L)^{\otimes l}$ on $U$ itself. Therefore we obtain a 
degree $d$ homogenous polynomial law
$$
\xi^* \mathcal{A}_l \to N(L)^{\otimes l}
$$
and therefore a morphism of $\mathcal{O}_S$-modules
$$
\Omega_l: \Gamma^d_{\mathcal{O}_T}(\xi^* \mathcal{A}_l) \to N(L)^{\otimes l}
$$
which we call the $l$-th \textit{Chow form} of $(Y, \mathcal{N})$.
It is easy to see that for any point $t \in T$ a local section 
$\phi$ of $\mathcal{A}_l$ which does not vanish at $t$ gives
a local section of $N(L)^{\otimes l}$ 
which does not vanish at $t$. Therefore $\Omega_l$ is a surjective
morphism of sheaves. Moreover, by multiplicativity of $N$
for a section $\phi_l$ of $\mathcal{A}_k$ and a section $\phi_m$ of
$\mathcal{A}_m$ we have equality
$$
\Omega_{m+l} (\phi_l \phi_m) = \Omega_{l}(\phi_l) \otimes \Omega_m(\phi_m)
$$
of local sections of $N(L)^{\otimes (m+l)}$. Therefore we obtain a 
surjective morphism of $\mathcal{O}_S$-algebras
$$
\xi^* \Big(\bigoplus_{l \geq 0} \Gamma^d_{\mathcal{O}_S}(\mathcal{A}_l) \Big)
\simeq 
\bigoplus_{l \geq 0} \Gamma^d_{\mathcal{O}_T}
(\xi^* \mathcal{A}_l) \to \bigoplus_{l \geq 0} N(L)^{\otimes l}
$$
hence an $S$-morphism 
$$
T \to Proj(\bigoplus_{l \geq 0} \Gamma^d_{\mathcal{O}_S}(\mathcal{A}_l))
$$
\begin{lemma} 
In the situation described above, $Proj(\bigoplus_{l \geq 0}
\Gamma^d_{\mathcal{O}_S}(\mathcal{A}_l)) \simeq \Gamma^d(X/S)$.
Moreover, if $l\geq r(d-1)$ then $\Gamma^d(X/S)$ is isomorphic to a closed subscheme of 
$\mathbb{P}(\Gamma^d(\mathcal{A}_l))$. When $S$ is a scheme over $\mathbb{Q}$, 
the assertion holds for any $l\geq 1$. 
\end{lemma}
\textit{Proof.} See Corollary 3.2.8 and Proposition 3.2.9 in \cite{R3}. 
$\square$
\begin{cor}
A family of cycles
$(Y, \mathcal{N})$ is
uniquely determined by its $l$-th Chow form
$\Omega_l: \Gamma^d_{\mathcal{O}_T}(\mathcal\xi^*{A}_l) \to N(L)^{\otimes l}$ where $l$ 
is given by the previous lemma.
\end{cor}

\bigskip
\noindent
\textit{Email:} vbaranov@math.uci.edu

\noindent
\textit{Address:} Department of Mathematics, UC Irvine, Irvine, CA 92697.

\end{document}